\documentclass[submission]{eptcs}
\providecommand{\event}{ACT 2021} 

\usepackage{underscore}
\usepackage{randtext}
\usepackage{setspace}

\usepackage[utf8]{inputenc}
\usepackage[T1]{fontenc}
\usepackage{hyperref}

\usepackage{verbatim}
\usepackage{graphicx}
\usepackage{rotating}

\usepackage{url}
\usepackage{hyperref}

\usepackage{pslatex}

\usepackage{datetime}

\usepackage{ebproof}

\usepackage{pgf,tikz,environ}
\usetikzlibrary{arrows}
\usetikzlibrary{cd}
\usetikzlibrary{decorations}

\usepackage{amsthm}
\usepackage[fleqn]{amsmath} 
\usepackage{amsfonts}
\usepackage{amssymb}
\usepackage{stix}

\theoremstyle{plain}
\newtheorem{theorem}{Theorem}[section]

\theoremstyle{definition}
\newtheorem{definition}[theorem]{Definition}
\newtheorem{example}[theorem]{Example}

\newcommand{\SetComp}[2]{\left\{ {#1}\:\middle|\:{#2} \right\}}

\definecolor{gold}{RGB}{255,215,0}
\definecolor{softBlack}{RGB}{45, 47, 49}
\definecolor{creamWhite}{RGB}{245,244,241}
\definecolor{softGray}{RGB}{220, 216, 214}
\definecolor{brick}{RGB}{232, 48, 48}

\usepackage{enumitem}

\title{Treewidth via Spined Categories (extended abstract)}
\author{Zoltan A. Kocsis
\institute{Commonwealth Scientific and Industrial Research Organisation\\ Eveleigh NSW, Australia}
\and
Benjamin Merlin Bumpus
\institute{
University of Glasgow\\
Scotland, UK}
}
\def\titlerunning{Treewidth via Spined Categories}
\def\authorrunning{Z. A. Kocsis \& B. M. Bumpus}
\begin{document}
\maketitle

\small

Treewidth is a well-known graph invariant with multiple interesting applications in combinatorics. On the practical side, many NP-complete problems are polynomial-time (sometimes even linear-time) solvable on graphs of bounded treewidth~\cite{cygan2015parameterized, flum2006parameterized}. On the theoretical side, treewidth played an essential role in the proof of the celebrated Robertson-Seymour graph minor theorem~\cite{RobertsonXX}. While defining treewidth-like invariants on graphs \cite{courcelle1993, HABIB201041, RobertsonX} and treewidth analogues on other sorts of combinatorial objects (incl. hypergraphs, digraphs~\cite{dtw, kreutzer2018}) has been a fruitful avenue of research, a direct, categorial description capturing multiple treewidth-like invariants is yet to emerge.

Here we report on our recent work on \textit{spined categories} \cite{mainArticle}: categories equipped with extra structure that permits the definition of a functorial analogue of treewidth, the \textit{triangulation functor}. The usual notion of treewidth is recovered as a special case, the triangulation functor of a spined category with graphs as objects and graph monomorphisms as arrows. The usual notion of treewidth for hypergraphs arises as the triangulation functor of a similar category of hypergraphs.

\section{Spined Categories}

Contrary to the usual convention in category-theoretic texts, we use the word \textit{graph} to refer to simple graphs (irreflexive, without loops or multiedges). We write $\mathrm{Grph}$ for the category that has graphs as objects and graph homomorphisms as arrows, and $\mathrm{Grph}_{m}$ for the category with the same objects, but monomorphisms as arrows.

\begin{definition}\label{def:spined-cat}
{ 
A \textit{spined category} consists of a category $\mathcal{C}$ equipped with the following additional structure:
\begin{itemize}\itemsep0em 
    \item a sequence $\Omega: \mathbb{N} \rightarrow \mathrm{ob}\:\mathcal{C}$ called the \textit{spine} of $\mathcal{C}$,
    \item an operation $\mathfrak{P}$ (called the \textit{proxy pushout}) that assigns to each diagram of the form \\
\begin{tikzcd}
G & \Omega_n \arrow[l, "g"'] \arrow[r, "h"] & H
\end{tikzcd}
    in $\mathcal{C}$ a distinguished cocone
\begin{tikzcd}
G \arrow[r, "{\mathfrak{P}(g,h)_g}"] & {\mathfrak{P}(g,h)} & H \arrow[l, "{\mathfrak{P}(g,h)_h}"'],
\end{tikzcd}
\end{itemize}
subject to the following two conditions:
\begin{enumerate}[label=\textbf{SC\arabic*}]\itemsep0em
    \item For every $X \in \mathrm{ob}\:\mathcal{C}$ there is $n \in \mathbb{N}$ such that $\mathcal{C}(X,\Omega_n) \neq \emptyset$.\label{property:SC1}
    \item Given any diagram of the form
\begin{tikzcd}
G' & G \arrow[l, "g'"'] & \Omega_n \arrow[l, "g"'] \arrow[r, "h"] & H \arrow[r, "h'"] & H'
\end{tikzcd}
    we can find a \emph{unique} morphism $(g',h'): \mathfrak{P}(g,h) \rightarrow \mathfrak{P}(g' \circ g, h' \circ h)$ making the following diagram commute:
\begin{center}
\begin{tikzcd}
\Omega_n \arrow[d, "h"'] \arrow[r, "g"]                            & G \arrow[r, "g'"] \arrow[d, "{\mathfrak{P}(g,h)_g}"] & G' \arrow[dd, "{\mathfrak{P}(g'\circ g,h'\circ h)_{g' \circ g}}"] \\
H \arrow[d, "h'"'] \arrow[r, "{\mathfrak{P}(g,h)_h}"']             & {\mathfrak{P}(g,h)} \arrow[rd, "{(g',h')}", dashed]  &                                                                   \\
H' \arrow[rr, "{\mathfrak{P}(g'\circ g,h'\circ h)_{h' \circ h}}"'] &                                                      & {\mathfrak{P}(g' \circ g,h' \circ h)}                            
\end{tikzcd}
\end{center}\label{property:SC2}
\end{enumerate}
}\end{definition}

Proxy pushouts capture an important property that the left-cancellative subcategory $\mathrm{Grph}_{m}$ "remembers" about the existence of pushouts in the category $\mathrm{Grph}$: we can equip the former with proxy pushouts by assigning to each diagram 
\begin{tikzcd}
G & \Omega_n \arrow[l, "g"', hook'] \arrow[r, "h", hook] & H
\end{tikzcd}
its pushout square in the latter. Moreover, a category $\mathcal{C}$ with all pushouts, when equipped with a sequence $\Omega: \mathbb{N} \rightarrow \mathrm{ob}\:\mathcal{C}$ satisfying \textbf{SC1}, always forms a spined category.

\begin{definition}\label{def:spinal-functor}
Consider spined categories $(\mathcal{C},\Omega^\mathcal{C}, \mathfrak{P}^\mathcal{C})$ and $(\mathcal{D},\Omega^\mathcal{D}, \mathfrak{P}^\mathcal{D})$. We call a functor $F: \mathcal{C} \rightarrow \mathcal{D}$ a \textit{spinal functor} if it
\begin{enumerate}[label=\textbf{SF\arabic*}]
    \item \textit{preserves the spine}, i.e. $F \circ \Omega^\mathcal{C} = \Omega^\mathcal{D}$, and \label{property:SF1}
    \item \textit{preserves proxy pushouts}, i.e. given a proxy pushout square
\begin{center}
\begin{tikzcd}
\Omega_n \arrow[d, "h"'] \arrow[r, "g"] & G \arrow[d, "{\mathfrak{P}^\mathcal{C}(g,h)_g}"] \\
H \arrow[r, "{\mathfrak{P}^\mathcal{C}(g,h)_h}"']   & {\mathfrak{P}^\mathcal{C}(g,h)}                 
\end{tikzcd}
\end{center}
    in the category $\mathcal{C}$, its $F$-image
\begin{center}
\begin{tikzcd}
\Omega_n \arrow[d, "Fh"'] \arrow[r, "Fg"]   & {F[G]} \arrow[d, "{F\mathfrak{P}^\mathcal{C}(g,h)_g}"] \\
{F[H]} \arrow[r, "{F\mathfrak{P}^\mathcal{C}(g,h)_h}"'] & {F[\mathfrak{P}^\mathcal{C}(g,h)]}
\end{tikzcd}
\end{center}
    forms a proxy pushout square in $\mathcal{D}$. One can state this equationally, by demanding that the equalities $F[\mathfrak{P}^\mathcal{C}(g,h)] = \mathfrak{P}^\mathcal{D}(Fg,Fh)$, $F\mathfrak{P}^\mathcal{C}(g,h)_g = \mathfrak{P}^\mathcal{D}(Fg,Fh)_{Fg}$ and $F\mathfrak{P}^\mathcal{C}(g,h)_h = \mathfrak{P}^\mathcal{D}(Fg,Fh)_{Fh}$ all hold.\label{property:SF2}
\end{enumerate}
\end{definition}

That is, a spinal functor between spined categories is a functor between the underlying categories that respects the spine and proxy pushout structure. As expected, the composition of two spinal functors is itself spinal.

Regard the poset $(\mathbb{N}, \leq)$ of natural numbers under the usual ordering as a category. This category has all pushouts. Equipping $(\mathbb{N}, \leq)$ with the spine $\Omega_n = n$ and suprema as proxy pushouts yields a simple example of a spined category, which we will denote $\mathrm{Nat}$.

\begin{definition}
An \textit{S-functor} on the spined category $C$ is a spinal functor defined on $C$ and valued in $\mathrm{Nat}$.
\end{definition}

Some spined categories do not have any S-functors defined on them: typically when some object $\Omega_n$ can be constructed as a proxy pushout using $\Omega_i$ for $i < n$. The interested reader is welcome to enumerate necessary/sufficient conditions for the existence of S-functors. In what follows, we side-step this issue by focusing our attention on the class of spined categories which have at least one S-functor defined on them. We call such categories \textit{measurable}.

\begin{example}
The category $\mathrm{Grph}_{m}$ (with proxy pushouts inherited from pushouts in $\mathrm{Grph}$, and spine $\Omega_n$ the complete\footnote{A graph where every pair of distinct vertices is connected by an edge.} graph on $n$ vertices) is measurable: it's easy to check that the map $\omega(G)$ which sends each $G$ to the size of its largest complete subgraph, constitutes an S-functor.  
\end{example}

\section{Triangulation Functor}

Our main result proves the existence of a distinguished S-functor, the \textit{triangulation functor} on each measurable spined category. Treewidth is recovered as the triangulation functor of the category $\mathrm{Grph}_{m}$, while hypergraph treewidth is recovered as the triangulation functor of a corresponding category $\mathrm{HGrph}_{m}$. For traditional graph-theoretic definitions of treewidth, we refer the reader to Encyclopedia~of~Algorithms~\cite{Bodlaender2016tw}: our \textit{pseudo-chordal} objects play a similar role to that of chordal\footnote{A graph where all cycles of $>3$ vertices have a \textit{chord}, i.e. an edge connecting non-adjacent vertices of the cycle.} graphs in the second characterisation presented there.

\begin{definition}\label{def:pseudo-chordal-object}
We call an object $X$ of a spined category $\mathcal{C}$ \emph{pseudo-chordal} if all S-functors assign the same value to $X$, i.e. for any two S-functors $F,G: \mathcal{C} \rightarrow \mathrm{Nat}$ we have $F[X] = G[X]$. We let $\mathrm{pc}\:\mathcal{C}$ denote the class of pseudo-chordal objects in the category $\mathcal{C}$.
\end{definition}

In the category $\mathrm{Grph}_{m}$ defined above, the class of pseudo-chordal objects forms a strict superset of the class of chordal graphs: while all chordal graphs are in fact pseudo-chordal objects, the converse fails.

\begin{theorem}[Main result]\label{thm:main}
Take a measurable spined category $\mathcal{C}$, equipped with some S-functor\\
$s: \mathcal{C} \rightarrow \mathrm{Nat}$. The map $\Delta: \mathcal{C} \rightarrow \mathrm{Nat}$ defined by the equation $\Delta[G] =  \min \SetComp{s(H)}{H \in \mathrm{pc}\:C, \mathcal{C}(G,H) \neq \emptyset}$
\begin{enumerate}
    \item is an S-functor;
    \item dominates all other S-functors, i.e. for any $X \in \mathrm{ob}\:\mathcal{C}$ and S-functor $F: \mathcal{C} \rightarrow \mathrm{Nat}$, $F[X] \leq \Delta[X]$.
\end{enumerate}
\end{theorem}

We call the functor $\Delta$ the triangulation functor of the category $\mathcal{C}$. It's clear that every measurable category has a unique triangulation functor.

\begin{theorem}
The triangulation functor of the category $\mathrm{Grph}_{m}$ coincides with treewidth.
\end{theorem}

Spined categories socialize well via spinal functors: in the talk, we will explain how one can obtain measurability (and non-measurability) results purely by constructing spinal functors, and present further examples, including a category of hypergraphs where the triangulation functor recovers the notion of hypergraph treewidth. Some previously unknown  tree-width-like invariants also emerge by collecting the relevant combinatorial objects into a spined category. Somewhat surprisingly, by putting mild computability conditions on the category, we can even obtain an algorithm which computes the value of the triangulation functor (although the generic algorithms obtained this way are impractically slow for computing the treewidth of all but the simplest graphs).

\begin{figure}[h!]
\begin{center}
\includegraphics[scale=0.2]{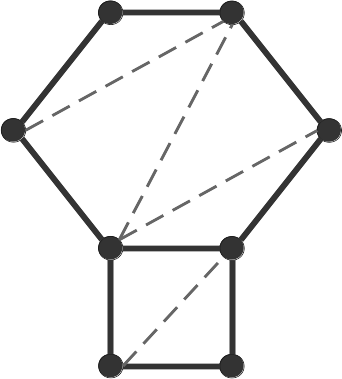}
\vspace*{-0.65cm}
\end{center}
\small
\caption{A graph and one of its chordal completions (dashed). The treewidth of a graph $G$ is the infimum of the sizes of the largest complete subgraphs contained in the chordal completions of $G$.}
\end{figure}

\vspace*{-0.8cm}

\bibliographystyle{eptcs}
\bibliography{biblio}

\begin{thebibliography}{10}
\providecommand{\bibitemdeclare}[2]{}
\providecommand{\surnamestart}{}
\providecommand{\surnameend}{}
\providecommand{\urlprefix}{Available at }
\providecommand{\url}[1]{\texttt{#1}}
\providecommand{\href}[2]{\texttt{#2}}
\providecommand{\urlalt}[2]{\href{#1}{#2}}
\providecommand{\doi}[1]{doi:\urlalt{http://dx.doi.org/#1}{#1}}
\providecommand{\bibinfo}[2]{#2}

\bibitemdeclare{inbook}{Bodlaender2016tw}
\bibitem{Bodlaender2016tw}
\bibinfo{author}{H.~L. \surnamestart Bodlaender\surnameend}
  (\bibinfo{year}{2016}): \emph{\bibinfo{title}{Treewidth of Graphs}}, pp.
  \bibinfo{pages}{2255--2257}.
\newblock \bibinfo{publisher}{Springer New York}, \bibinfo{address}{New York,
  NY}, \doi{https://doi.org/10.1007/978-1-4939-2864-4_431}.

\bibitemdeclare{article}{mainArticle}
\bibitem{mainArticle}
\bibinfo{author}{B.~M. \surnamestart {Bumpus}\surnameend} \&
  \bibinfo{author}{Z.~A. \surnamestart {Kocsis}\surnameend}
  (\bibinfo{year}{2021}): \emph{\bibinfo{title}{{Spined categories:
  generalizing tree-width beyond graphs}}}.
\newblock {\sl \bibinfo{journal}{arXiv
  e-prints}}:\bibinfo{eid}{arXiv:2104.01841}.

\bibitemdeclare{article}{courcelle1993}
\bibitem{courcelle1993}
\bibinfo{author}{B.~\surnamestart Courcelle\surnameend},
  \bibinfo{author}{J.~\surnamestart Engelfriet\surnameend} \&
  \bibinfo{author}{G.~\surnamestart Rozenberg\surnameend}
  (\bibinfo{year}{1993}): \emph{\bibinfo{title}{Handle-rewriting hypergraph
  grammars}}.
\newblock {\sl \bibinfo{journal}{Journal of Computer and System Sciences}}
  \bibinfo{volume}{46}(\bibinfo{number}{2}), pp. \bibinfo{pages}{218--270},
  \doi{https://doi.org/10.1016/0022-0000(93)90004-G}.

\bibitemdeclare{book}{cygan2015parameterized}
\bibitem{cygan2015parameterized}
\bibinfo{author}{M.~\surnamestart Cygan\surnameend}, \bibinfo{author}{F.~V.
  \surnamestart Fomin\surnameend}, \bibinfo{author}{{\L}.~\surnamestart
  Kowalik\surnameend}, \bibinfo{author}{D.~\surnamestart
  Lokshtanov\surnameend}, \bibinfo{author}{D.~\surnamestart Marx\surnameend},
  \bibinfo{author}{M.~\surnamestart Pilipczuk\surnameend},
  \bibinfo{author}{M.~\surnamestart Pilipczuk\surnameend} \&
  \bibinfo{author}{S.~\surnamestart Saurabh\surnameend} (\bibinfo{year}{2015}):
  \emph{\bibinfo{title}{Parameterized algorithms}}.
\newblock \bibinfo{publisher}{Springer},
  \doi{https://doi.org/10.1007/978-3-319-21275-3}.

\bibitemdeclare{article}{flum2006parameterized}
\bibitem{flum2006parameterized}
\bibinfo{author}{J.~\surnamestart Flum\surnameend} \&
  \bibinfo{author}{M.~\surnamestart Grohe\surnameend} (\bibinfo{year}{2006}):
  \emph{\bibinfo{title}{Parameterized Complexity Theory. 2006}}.
\newblock {\sl \bibinfo{journal}{Texts Theoret. Comput. Sci. EATCS Ser}},
  \doi{https://doi.org/10.1007/3-540-29953-X}.

\bibitemdeclare{article}{HABIB201041}
\bibitem{HABIB201041}
\bibinfo{author}{M.~\surnamestart Habib\surnameend} \&
  \bibinfo{author}{C.~\surnamestart Paul\surnameend} (\bibinfo{year}{2010}):
  \emph{\bibinfo{title}{A survey of the algorithmic aspects of modular
  decomposition}}.
\newblock {\sl \bibinfo{journal}{Computer Science Review}}
  \bibinfo{volume}{4}(\bibinfo{number}{1}), pp. \bibinfo{pages}{41 -- 59},
  \doi{https://doi.org/10.1016/j.cosrev.2010.01.001}.

\bibitemdeclare{article}{dtw}
\bibitem{dtw}
\bibinfo{author}{T.~\surnamestart Johnson\surnameend},
  \bibinfo{author}{N.~\surnamestart Robertson\surnameend},
  \bibinfo{author}{P.~D. \surnamestart Seymour\surnameend} \&
  \bibinfo{author}{R.~\surnamestart Thomas\surnameend} (\bibinfo{year}{2001}):
  \emph{\bibinfo{title}{Directed tree-width}}.
\newblock {\sl \bibinfo{journal}{Journal of Combinatorial Theory. Series B}}
  \bibinfo{volume}{82}(\bibinfo{number}{1}), pp. \bibinfo{pages}{138--154},
  \doi{https://doi.org/10.1006/jctb.2000.2031}.

\bibitemdeclare{incollection}{kreutzer2018}
\bibitem{kreutzer2018}
\bibinfo{author}{S.~\surnamestart Kreutzer\surnameend} \&
  \bibinfo{author}{O.-j. \surnamestart Kwon\surnameend} (\bibinfo{year}{2018}):
  \emph{\bibinfo{title}{Digraphs of Bounded Width}}.
\newblock In: {\sl \bibinfo{booktitle}{Classes of Directed Graphs}},
  \bibinfo{publisher}{Springer}, pp. \bibinfo{pages}{405--466},
  \doi{https://doi.org/10.1007/978-3-319-71840-8_9}.

\bibitemdeclare{article}{RobertsonX}
\bibitem{RobertsonX}
\bibinfo{author}{N.~\surnamestart Robertson\surnameend} \&
  \bibinfo{author}{P.~D. \surnamestart Seymour\surnameend}
  (\bibinfo{year}{1991}): \emph{\bibinfo{title}{Graph minors X. Obstructions to
  tree-decomposition}}.
\newblock {\sl \bibinfo{journal}{Journal of Combinatorial Theory, Series B}}
  \bibinfo{volume}{52}(\bibinfo{number}{2}), pp. \bibinfo{pages}{153--190},
  \doi{https://doi.org/10.1016/0095-8956(91)90061-N}.

\bibitemdeclare{article}{RobertsonXX}
\bibitem{RobertsonXX}
\bibinfo{author}{N.~\surnamestart Robertson\surnameend} \&
  \bibinfo{author}{P.D. \surnamestart Seymour\surnameend}
  (\bibinfo{year}{2004}): \emph{\bibinfo{title}{Graph Minors. XX. Wagner's
  conjecture}}.
\newblock {\sl \bibinfo{journal}{Journal of Combinatorial Theory, Series B}}
  \bibinfo{volume}{92}(\bibinfo{number}{2}), pp. \bibinfo{pages}{325 -- 357},
  \doi{https://doi.org/10.1016/j.jctb.2004.08.001}.

\end{thebibliography}
\end{document}